\documentclass{amsart}

\usepackage{amsmath,amssymb,amsthm}
\usepackage{graphicx,tikz}

\newtheorem*{thm}{Theorem}

\theoremstyle{definition}

\theoremstyle{remark}

\begin{document}

\title[]{\textsc{Max-Cut} via Kuramoto-type Oscillators}
\subjclass[2010]{} 
\thanks{S.S. is supported by the NSF (DMS-1763179) and the Alfred P. Sloan Foundation.}

\author[]{Stefan Steinerberger}
\address{Department of Mathematics, University of Washington, Seattle}
\email{steinerb@uw.edu}

\begin{abstract} We consider the \textsc{Max-Cut} problem. Let $G = (V,E)$ be a graph with adjacency matrix $(a_{ij})_{i,j=1}^{n}$. Burer, Monteiro \& Zhang proposed to find, for $n$ angles $\left\{\theta_1, \theta_2, \dots, \theta_n\right\} \subset [0, 2\pi]$, minima of the energy
$$ f(\theta_1, \dots, \theta_n) = \sum_{i,j=1}^{n} a_{ij} \cos{(\theta_i - \theta_j)}$$
because configurations achieving a global minimum leads to a partition of size $0.878\cdot\textsc{Max-Cut}(G)$. This approach is known to be computationally viable and leads to very good results in practice. We prove that by replacing $\cos{(\theta_i - \theta_j)}$ with an explicit function $g_{\varepsilon}(\theta_i - \theta_j)$ global minima of this new functional lead to a $(1-\varepsilon)\cdot\textsc{Max-Cut}(G)$. This suggests some interesting algorithms that perform well. It also shows that the problem of finding approximate global minima of energy functionals of this type is NP-hard in general. 
 \end{abstract}

\maketitle

\section{Introduction}
\subsection{Max-Cut} We consider a classical problem: given a graph $G=(V,E)$, what is the best decomposition of its vertices into two sets such that the number
of edges between the two sets is maximal? \textsc{Max-Cut} is known to be NP-hard.

\begin{center}
\begin{figure}[h!]
\begin{tikzpicture}[scale=0.8]
\filldraw (0,0) circle (0.06cm);
\filldraw (-0.5,1) circle (0.06cm);
\filldraw (0.5,2) circle (0.06cm);
\filldraw (0.5,1) circle (0.06cm);
\filldraw (3,0) circle (0.06cm);
\filldraw (3.2-0.5,1.1) circle (0.06cm);
\filldraw (3.7,1.8) circle (0.06cm);
\filldraw (3.5,1.2) circle (0.06cm);
\draw [dashed] (0,1) ellipse (1cm and 1.6cm);
\draw [dashed] (3.2,1) ellipse (1cm and 1.6cm);
\draw [thick] (0.8,2) -- (2.3,1.7);
\draw [thick] (1,1) -- (2.2,1.5);
\draw [thick] (0.8,0) -- (2.2,0.5);
\draw [thick] (1,0.5) -- (2.2,0.3);
\draw [thick] (1,1) -- (2.2,1.1);
\end{tikzpicture}
\caption{Decomposing Vertices of a Graph into two sets so that many edges run between them.}
\end{figure}
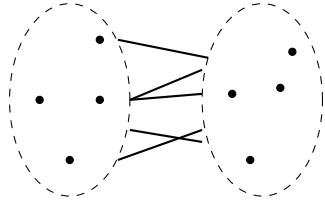
\end{center}

It is easy to see that by picking the subsets uniformly at random, we will get, in expectation, a partition such that at least $|E|/2$ edges run between them. This shows that it is easy to find a partition such that $0.5 \cdot \textsc{Max-Cut}(G)$ edges run between them. In a seminal paper, Goemans \& Williamson \cite{goemans} constructed a $0.878-$approximation algorithm for \textsc{Max-Cut}, where the constant is given by
$$ 0.878 \dots = \frac{2}{\pi} \min_{0 \leq \theta \leq \pi} \frac{\theta}{1- \cos{(\theta)}}.$$
The algorithm uses semi-definite programming and randomized rounding (randomized rounding is explained in greater detail below). If the Unique Games Conjecture is true, this is the best possible approximation ratio for \textsc{Max-Cut} \cite{khot} that can be computed in polynomial time. It is known unconditionally that approximating \textsc{Max-Cut} by any factor better than $16/17 \sim 0.941$ is also NP-hard \cite{bell, hastad, trev}.

\subsection{Kuramoto Oscillators.} Kuramoto Oscillators refers to a broad class of problems where we are given $n$ particles $\theta_1, \dots, \theta_n \in \mathbb{S}^1 \cong [0, 2\pi]$.  These particles, which often depend on time, are assumed to be coupled in some nontrivial way (we refer to the surveys \cite{review1, review2}). A particularly nice setting is to define the energy
$$ f(\theta_1, \dots, \theta_n) = \sum_{i,j=1}^{n} a_{ij} \cos{(\theta_i - \theta_j)},$$
where $a_{ij} \in \left\{0,1\right\}$ is the entry of an adjacency matrix of a graph whose structure models the dependency between the particles. The particles are then assumed to move along the gradient of the energy, the overarching question is whether the underlying graph structure forces some type of universal behavior on the particles. The energy landscape of this particular energy, for example, is quite intricate. 
Taylor \cite{taylor} proved that if each vertex is connected to at least $\mu(n-1)$ vertices for $\mu \geq 0.9395$, then $f(\theta_1, \dots, \theta_n)$ does not have local maxima that are not also global. This was then improved by Ling, Xu \& Bandeira \cite{ling} to $\mu \geq 0.7929$ and J. Lu and the author \cite{jianfeng} to $\mu \geq 0.7889$. Townsend, Stillman \& Strogatz \cite{townsend} suggest that the critical value could be $\mu_c = 0.75$ --  they also identify networks with $\mu = 0.75$ having interesting spectral properties. These results mirror consideration by Burer, Monteiro \& Zhang \cite{burer3} and are partially inspired by those.

\subsection{The Approach of Burer, Monteiro \& Zhang} Burer, Monteiro \& Zhang \cite{burer3} proposed a particular rank-two relaxation of the Goemans-Williamson approach \cite{goemans}. We recall that Goemans-Williamson suggested to relax 
$$ 2 \cdot |E| - 4 \cdot \textsc{MaxCut}(G) =  \min_{x_i \in \left\{-1,1\right\}} \sum_{i,j=1}^{n} a_{ij} x_i x_j$$
by replacing the $x_i \in \left\{-1,1\right\}$ with unit vectors $v_i \in \mathbb{R}^n$ and $x_i x_j$ with $\left\langle v_i, v_j\right\rangle$. This is clearly a more general problem but one that is amenable to being solved with SDP methods in polynomial time. In the last step, they perform a randomized rounding step and prove that this leads to a $0.878 \cdot \textsc{Max-Cut}$ approximation. Burer, Monteiro \& Zhang \cite{burer3} suggest that it might be possible to bypass the SDP step by arguing directly on the relaxed problem in $\mathbb{R}^2$. Parametrizing unit vectors in $\mathbb{R}^2$ by an angle $\theta \in \mathbb{S}^1$, we see that
$$ \left\langle v_{\theta_i}, v_{\theta_j} \right\rangle = \cos{(\theta_i - \theta_j)}$$
and this leads to the notion of energy $f:(\mathbb{S}^1)^n \cong [0,2\pi]^n \rightarrow \mathbb{R}$
 $$ f(\theta_1, \dots, \theta_n) = \sum_{i,j=1}^{n} a_{ij} \cos{(\theta_i - \theta_j)}.$$
  Burer, Monteiro \& Zhang propose to minimize this energy instead and then use the same randomized rounding step as in the Goemans-Williamson approach.
The success of this particular relaxation will depend on two competing factors.
\begin{itemize}
\item \textbf{Upside.} There is no longer any need for solving a semi-definite program (which becomes computationally expensive when $n$ is large), one simply has to find a configuration of particles $\left\{\theta_1, \dots, \theta_n \right\} \subset [0,2\pi]$ for which the energy is small.  Moreover, there is no need to find the global minimum (but: the smaller the energy, the better the configuration). 
\item \textbf{Downside.} Without the SDP, there is no particular insight into how one would start looking for a good configuration  $\left\{\theta_1, \dots, \theta_n \right\} \subset [0,2\pi]$.  Gradient descent methods are at the mercy of the energy landscape, it might potentially be hard to find a configuration for which the energy is small.
\end{itemize}

In practice, the (hypothetical) downside does not seem to cause any difficulties, the Burer, Monteiro \& Zhang (BUR02) approach is known to work very well. Indeed, in an extensive 2018 comparison, Dunning, Gupta \& Silberholz \cite{dunning} compared 37 different heuristics over 3296 problem instances concluding:
``The best overall heuristic on the expanded instance library with
respect to the performance of its mean solution across
the five replicates for each instance was max-cut heuristic BUR02, which not only matched the best performance on 22.9\% of instances but also had strictly better performance than any other heuristic on 16.2\% of instances and a mean deviation of only 0.3\%.'' It is not entirely understood why this relaxation works so well and this is being actively studied, we refer to  Boumal, Voroninski, Bandeira \cite{boumal, boumal2}, Ling \cite{ling0} and Ling, Xu \& Bandeira \cite{ling}.

\section{The Result}
\subsection{Main Idea.} Our main idea is the following: the cosine arises naturally when considering the inner product between two vectors since
$$ \left\langle v_{\theta_i}, v_{\theta_j} \right\rangle = \cos{(\theta_i - \theta_j)}.$$
However, since we are not actually using any type of SDP approach, we do not really have to use the cosine. Maybe there are other functions $g:\mathbb{S}^1 \rightarrow [-1,1]$ that are as good or possibly even better? The only constraint is that the randomized rounding step, when applied to a minimal energy configuration, should work well.
We will consider more general notions of a Kuramoto-type energy of the form
 $$ f(\theta_1, \dots, \theta_n) = \sum_{i,j=1}^{n} a_{ij} \cdot g(\theta_i - \theta_j),$$
where $g: \mathbb{S}^1 \cong [0,2\pi] \rightarrow \mathbb{R}$ is assumed to 
\begin{enumerate}
\item be differentiable everywhere,
\item be symmetric in the sense of $g(x) = g(-x)$ and
\item to assume its maximum in $g(0) = 1$ and its minimum in $g(\pi) = -1$.
\end{enumerate}

Minimizing such a Kuramoto-type energy will encourage that any two vertices $v_1, v_2$ that are connected by an edge $(v_1, v_2) \in E$ are moved to antipodal points on the circle. If the underlying graph is bipartite, this will indeed be the unique minimal energy configuration. For more general graphs, this is not so simple and one would expect a minimal energy configuration to depend on the graph.
We want that minimal energy configurations are well-behaved with respect to randomized rounding. For any given set of angles $\left\{\theta_1, \dots, \theta_n\right\} \subset \mathbb{S}^1$, the randomized rounding procedure results in an assignment of points into two sets as follows  (see Fig. 2): pick a random line going through the origin which splits the sets into the two groups induced by two half-spaces and use those sets as a partition.

 \begin{center}
\begin{figure}[h!]
\begin{tikzpicture}[scale=0.9]
\draw [thick] (0,0) circle (2cm);
\filldraw (0,0) circle (0.04cm);
\filldraw (2,0) circle (0.06cm);
\filldraw (1.9,0.6) circle (0.06cm);
\filldraw (1.92,-0.5) circle (0.06cm);
\filldraw (-1.82,0.8) circle (0.06cm);
\filldraw (-1.92,-0.5) circle (0.06cm);
\filldraw (-2,0) circle (0.06cm);
\filldraw (1.4142, 1.4142) circle (0.06cm);
\node at (2.3, 0) {$\theta_4$};
\node at (2.2, -0.6) {$\theta_6$};
\node at (-2.2, -0.6) {$\theta_{2}$};
\node at (-2.3, 0) {$\theta_{5}$};
\node at (1.75, 1.41) {$\theta_{7}$};
\node at (2.3, 0.55) {$\theta_{3}$};
\node at (-2.2, 0.9) {$\theta_{1}$};
\draw[dashed] (-1.8, -2.2) -- (1.8, 2.2);
\draw (-2,0) ellipse (1cm and 1.5cm);
\draw (1.8,0.2) ellipse (1.2cm and 1.5cm);
\end{tikzpicture}
\caption{Randomized Rounding: for given $\left\{\theta_1, \dots, \theta_n\right\} \subset \mathbb{S}^1$, we can pick a random line through the origin and the partition the vertices of the Graph according to the two half-spaces.}
\end{figure}
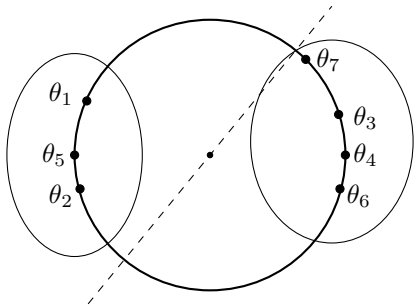
\end{center}
\vspace{-10pt}
We can analyze the expected behavior of randomized rounding completely in terms of the energy $f(\theta_1, \dots, \theta_n)$. For the minimal energy configuration of Kuramoto-type energies of this flavor, we obtain the following approximation result.

\begin{thm} Let $g:\mathbb{S}^1 \rightarrow \mathbb{R}$ be an admissible function and let $\left\{\theta_1, \theta_2, \dots, \theta_n \right\} \subset \mathbb{S}^1$ be a minimal energy configuration of the associated energy.
 Then the expected number of edges for a randomized rounding partition satisfies
$$ \mathbb{E}_{} ~\emph{edges} \geq \left(  \min_{0 \leq x \leq \pi}   \frac{2}{\pi}  \frac{x}{1 -g(x)}  \right) \cdot \emph{\textsc{Max-Cut}}(G).$$
\end{thm}

If $g(x) = \cos{(x)}$, we recover the classical $0.878 \cdot \textsc{Max-Cut}(G)$ result. However, for more general $g(x)$, the constant can be arbitrarily close to 1. We also show that the size of the energy functional has immediate implications for the quality of the randomized rounding step by proving the inequality
$$ \mathbb{E}_{} ~\mbox{edges} \geq  \left[\left(  \min_{0 \leq x \leq \pi}   \frac{2}{\pi}  \frac{x}{1 -g(x)}  \right)  \right] \cdot \left(\frac{|E|}{2} - \frac{1}{4} f(\theta_1, \dots, \theta_n)\right).$$
Thus, as in BUR02 \cite{burer3}, we do not necessarily need to find a global minimum, it suffices to find configurations with small energy (and the smaller, the better).
\begin{center}
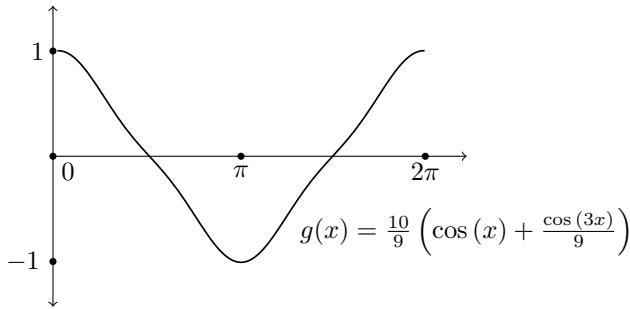
\begin{figure}[h!]
\begin{tikzpicture}
\node at (0,0) {\includegraphics[width=0.4\textwidth]{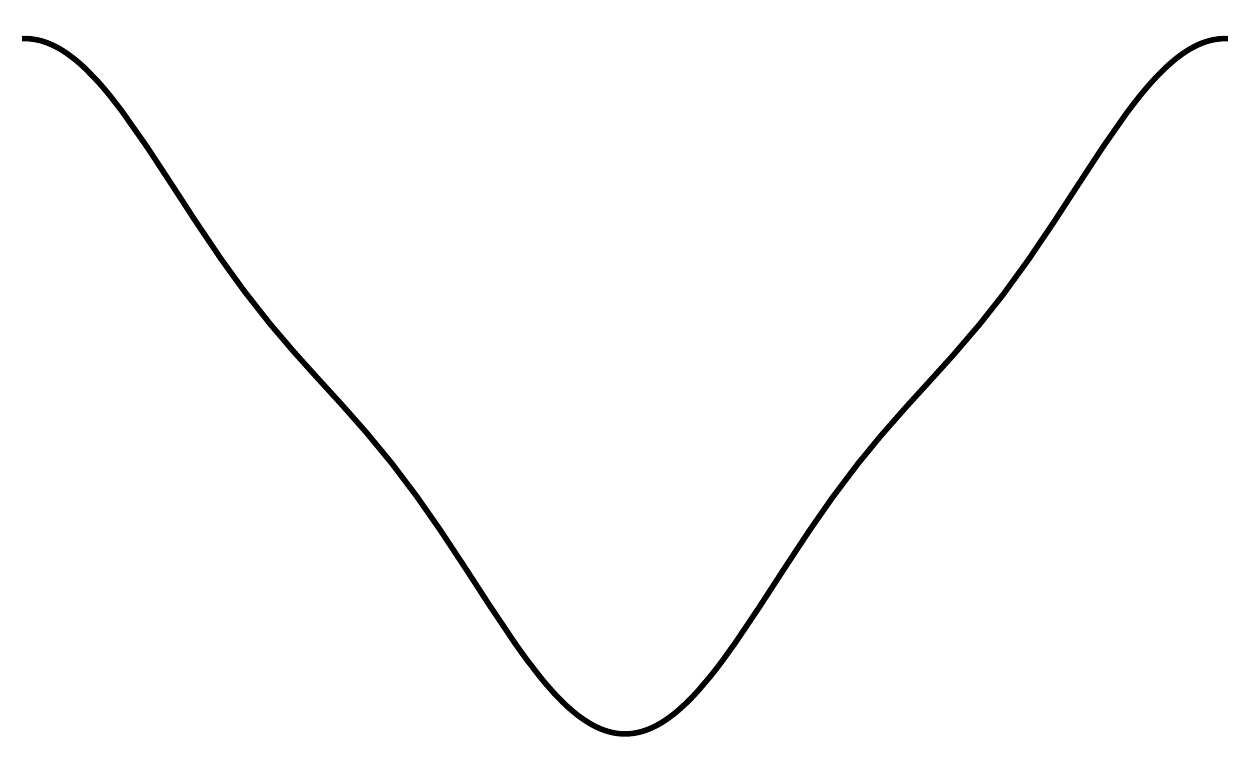}};
\draw [->] (-2.5,0) -- (3,0);
\draw [<->] (-2.5,-2) -- (-2.5,2);
\filldraw (-2.5, 0) circle (0.04cm);
\node at (-2.3, -0.2) {$0$};
\filldraw (-2.5, 1.4) circle (0.04cm);
\node at (-2.7, 1.4) {$1$};
\filldraw (-2.5, -1.4) circle (0.04cm);
\node at (-2.9, -1.4) {$-1$};
\filldraw (2.45, 0) circle (0.04cm);
\node at (2.45, -0.2) {$2\pi$};
\filldraw (0, 0) circle (0.04cm);
\node at (0, -0.2) {$\pi$};
\node at (3,-1) {$g(x) = \frac{10}{9} \left( \cos{(x)} + \frac{\cos{(3x)}}{9} \right)$};
\end{tikzpicture}
\vspace{-10pt}
\caption{Minima of this energy give $0.93517 \cdot \textsc{Max-Cut}$.}
\end{figure}
\end{center}

Our result is especially interesting when $\varepsilon > 0$ is small. Indeed, it is known \cite{bell, hastad, trev} that as soon as $\varepsilon < 1/17$ any $(1-\varepsilon)-$approximation of \textsc{Max-Cut} is necessarily NP-hard. In fact, if the unique games conjecture \cite{khot} is true, then the Goemans-Williamson approximation ratio of $0.878 \cdot \textsc{Max-Cut}$ is the best that one can do in polynomial time. This has an interesting consequence for the energy landscape of the energy $f(\theta_1, \dots, \theta_n)$ for functions $g$ for which the constant is  $> 16/17$: it must then, in general, be NP-hard to find a configuration $\left\{\theta_1, \dots, \theta_n\right\} \subset \mathbb{S}^1$ with energy close to the global minimum. We believe this to be an interesting statement about a large class of Kuramoto-type energy functionals. 

\subsection{Related results.} We are not aware of any results of this type. Most closely related in spirit is perhaps the idea of using oscillators to solve problems of this type \cite{chou, mallick, wang2}. Wang \& Roychowdhury \cite{wang}, for example, consider systems of coupled self-sustaining nonlinear oscillators. The main idea is that these are governed by a Lyapunov function that is closely related to the Ising Hamiltonian of the coupling graph which allows for approximations to \textsc{Max-Cut}.

\subsection{Examples.} We start with a completely explicit example and take an Erd\H{o}s-Renyi random graph $G(500, 0.01)$. The graph has $|V| = 500$ vertices and $|E| = 1549$ edges. Several runs of Goemans-Williamson (GW) show $\textsc{Max-Cut}(G) \geq 1176$. We start by minimizing the BUR02 energy
$$ f(\theta_1, \dots, \theta_n) = \sum_{i,j=1}^{n} a_{ij} \cos{(\theta_i - \theta_j)}$$
using a random initialization for the angles $\theta_1, \dots, \theta_n$ and standard gradient descent. The result is shown in Fig. 4. We see that the points seem to be  distributed all over the circle and we get somewhat nice uniform control: the arising cut is never too small and for certain angles clearly improves on the GW method.

 \begin{center}
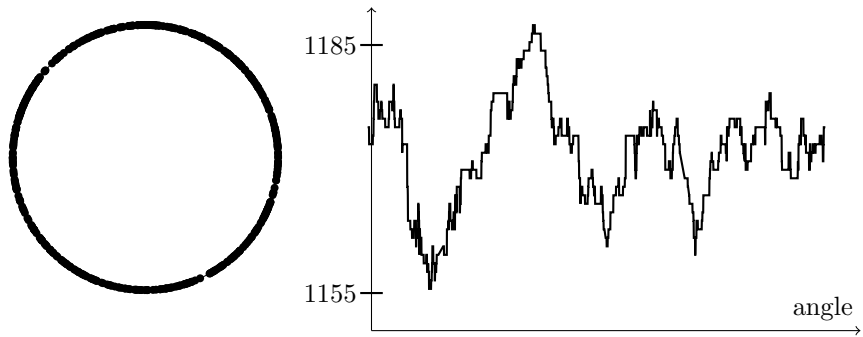
\begin{figure}[h!]
\begin{tikzpicture}
\node at (-1,0) {\includegraphics[width=0.3\textwidth]{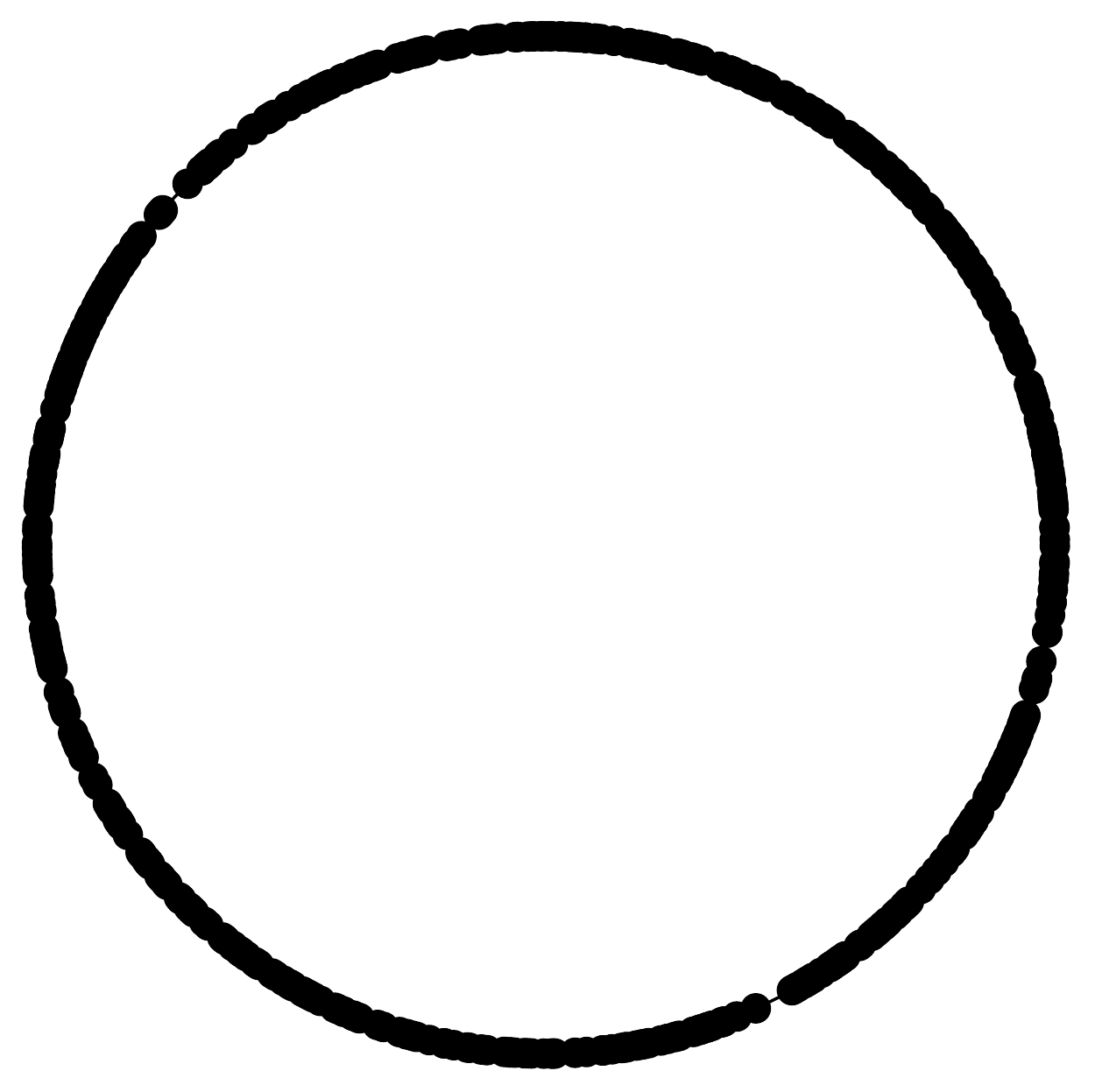}};
\node at (5,0) {\includegraphics[width=0.5\textwidth]{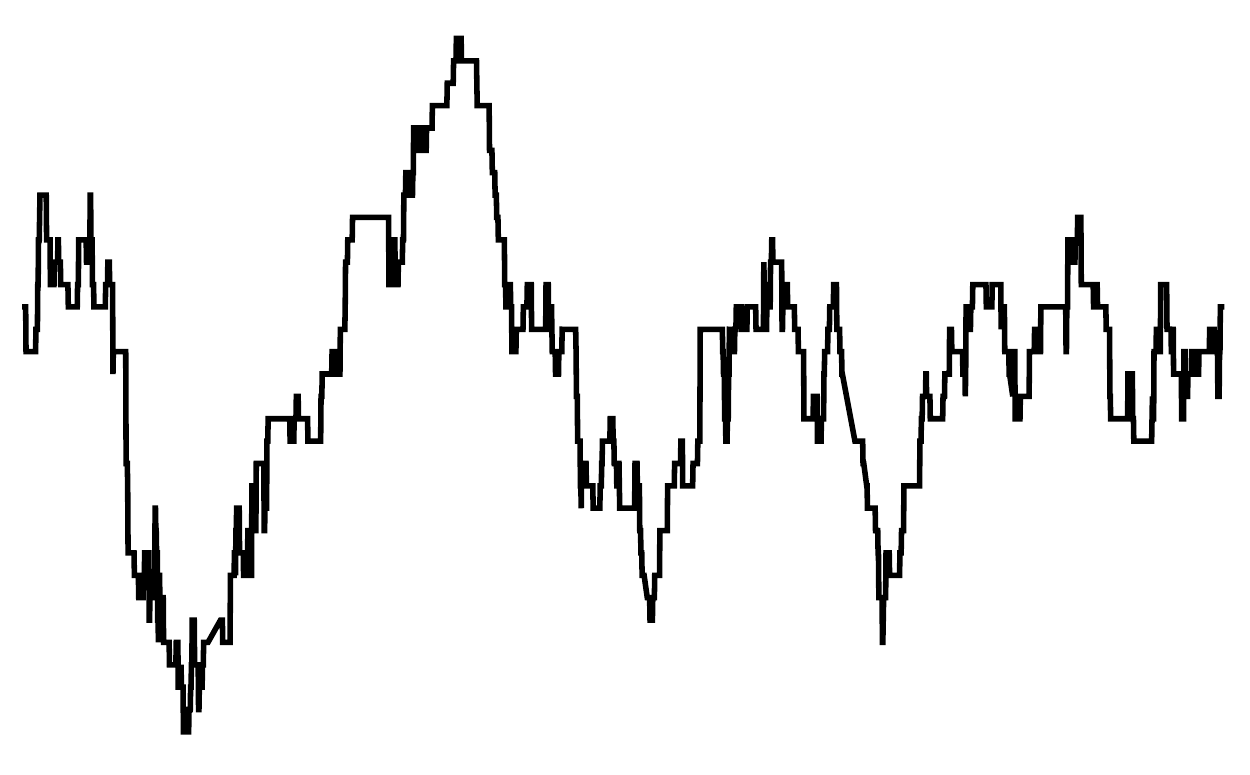}};
\draw [->] (2,-2.3) -- (8.5,-2.3);
\draw [->] (2,-2.3) -- (2,2);
\node at (8, -2) {$\mbox{angle}$};
\draw [thick] (1.85, 1.5) -- (2.15, 1.5);
\draw [thick] (1.85, -1.8) -- (2.15, -1.8);
\node at (1.45, 1.5) {1185};
\node at (1.45, -1.8) {1155};
\end{tikzpicture}
\caption{The distribution of points and the size of the cut obtained as a function of the angle of random line.}
\end{figure}
\end{center}

The question is now whether this can be improved by picking a function different from the cosine. There is a theoretical criterion (coming from Theorem 1) on how this function should look like
in the sense that
$$  \min_{0 \leq x \leq \pi}   \frac{2}{\pi}  \frac{x}{1 -g(x)}   \qquad \mbox{should be close to}~1.$$
There are many such functions -- a better understanding of which function $g(x)$ to choose would be interesting (see \S 2.5). We will use (throughout the paper)
$$ g(x) = \frac{99225}{117469}\left( \cos{(x)} + \frac{\cos{(3x)}}{9}  + \frac{\cos{(5x)}}{25} + \frac{\cos{(7x)}}{49} + \frac{\cos{(9x)}}{81} \right)$$
which comes from the Fourier series (normalized to $g(0) = 1 = - g(\pi)$) of  (see \S 2.5)
$$ 1 - \frac{2}{\pi} \cdot d_{\mathbb{S}^1}(0, x) = \begin{cases} 1 - \frac{2 x}{\pi} \qquad &\mbox{if}~0 \leq x \leq \pi \\
 1 - \frac{2 (2\pi -x)}{\pi} \qquad &\mbox{if}~\pi \leq x \leq 2\pi, \end{cases}$$
 where $d_{\mathbb{S}^1}(\cdot, \cdot)$ is the shortest distance on $\mathbb{S}^1$ (always less than $\pi$).
In particular, 
$$  \min_{0 \leq x \leq \pi}   \frac{2}{\pi}  \frac{x}{1 -g(x)} = 0.973$$
and every global minimum of this energy gives rise to a $0.973 \cdot \textsc{Max-Cut}$ approximation.
We run gradient descent (using the previously obtained final configuration of angles from the BUR02 method as initial set) and arrive at a nice result: the best cut has an additional 15 edges and all the cuts are uniformly closer to the maximum. Moreover, there is an additional `crystallization' of the points, hard to see in the picture, which are more structured (see \S 2.4).

 \begin{center}
\begin{figure}[h!]
\begin{tikzpicture}
\node at (-1,0) {\includegraphics[width=0.3\textwidth]{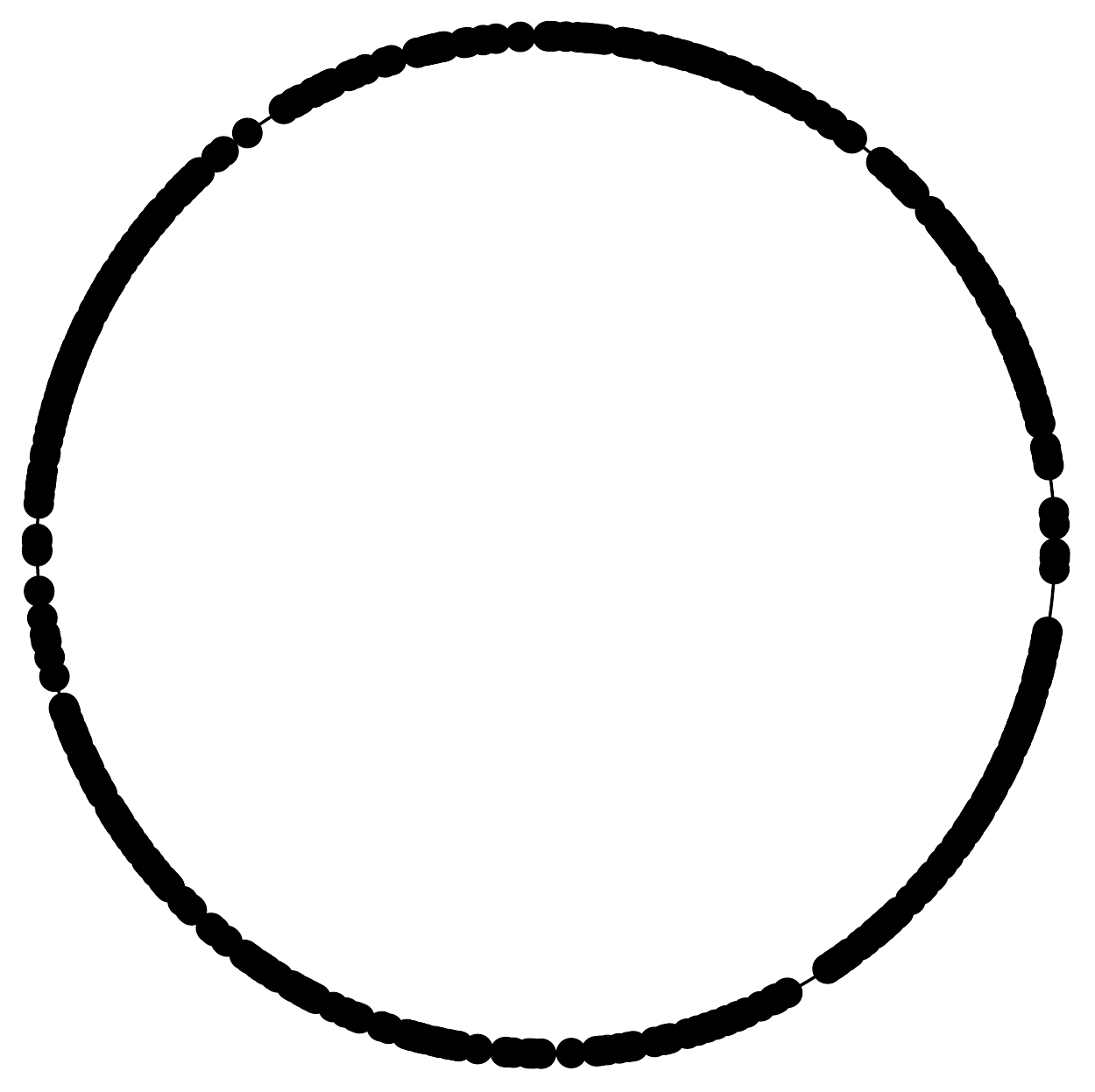}};
\node at (5,0) {\includegraphics[width=0.5\textwidth]{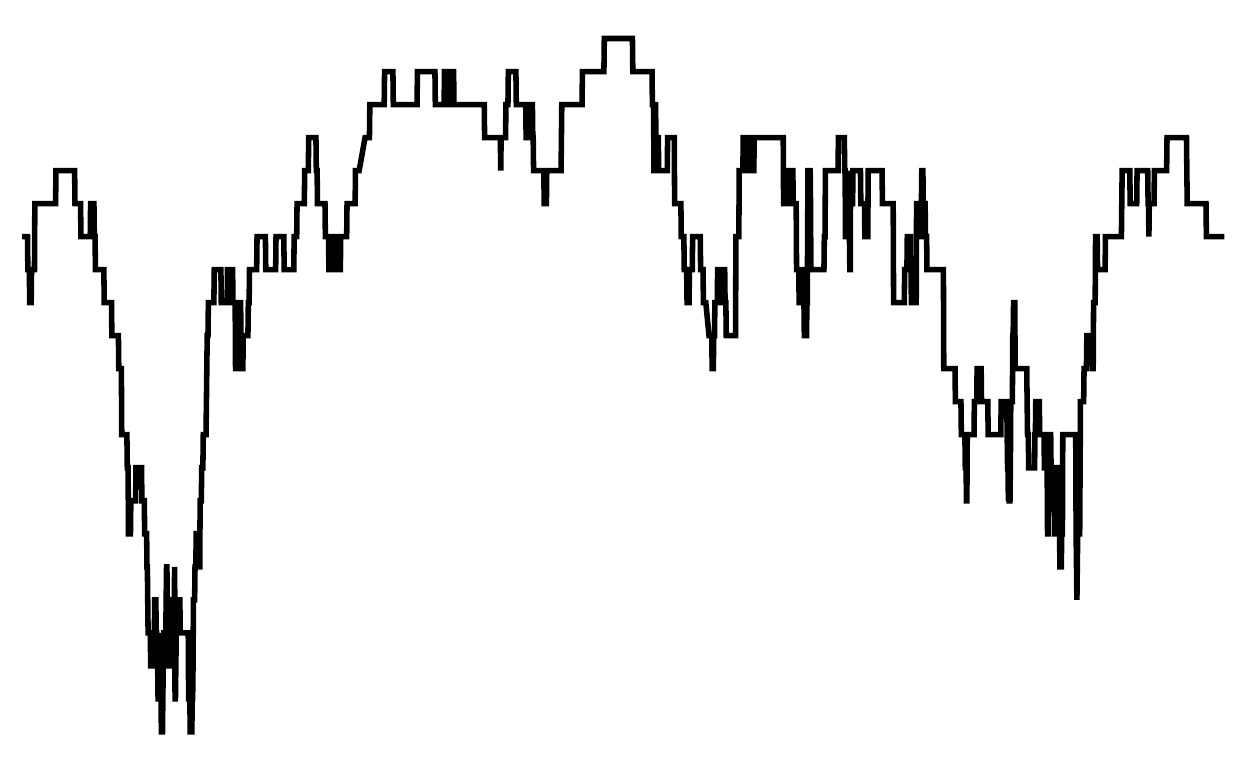}};
\draw [->] (2,-2.3) -- (8.5,-2.3);
\draw [->] (2,-2.3) -- (2,2);
\node at (8, -2) {$\mbox{angle}$};
\draw [thick] (1.85, 1.5) -- (2.15, 1.5);
\draw [thick] (1.85, -1.8) -- (2.15, -1.8);
\node at (1.45, 1.5) {1200};
\node at (1.45, -1.8) {1180};
\end{tikzpicture}
\caption{The distribution of points and the size of the cut obtained as a function of the angle of random line.}
\end{figure}
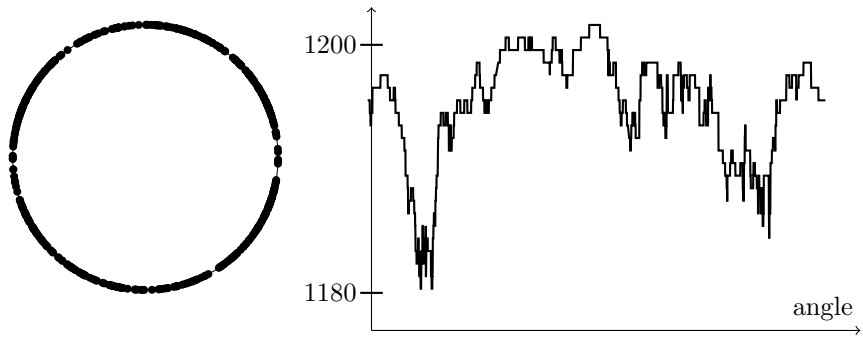
\label{fig:mont}
\end{center}

\vspace{-10pt}

This example appears to be quite typical for Erd\H{o}s-Renyi random graphs. Typically both the quality of the largest cut as well as the expected size of a random cut increases (our proof suggests why the expected size would increase). 

\begin{center}
\begin{table}[h!]
\begin{tabular}{ l c c c c r }
  Graph & $|V|$ & $ |E|$ & GW & BUR02 &  Our Method\\
  \hline
Mesner Graph $M_{22}$ & 77 & 616 & 400 & 420 & 420 \\
 Livingstone Graph & 266 & 1463 & 955 & 981 & 991 \\
 Berlekamp-Van Lint-Seidel & 243 & 2673 & 1572 & 1590 & 1606 \\
 Cameron Graph & 231 & 3465 & 1870 & 1884 & 1896 \\
\end{tabular}
\vspace{10pt}
\caption{Lower bounds on \textsc{Max-Cut} obtained by three methods.}
\end{table}
\end{center}
One could wonder whether these are artifacts coming from the randomness of the Erd\H{o}s-Renyi graphs. We decided to compare performance on some structured graphs for which we were unable to find the value of \textsc{Max-Cut} in the literature. Several runs of each method leads to the bounds on \textsc{Max-Cut} in Table 1.

\subsection{The Crystallization Phenomenon.} Minimal energy configurations of our functional tend to be somewhat structured. We start with an example.  The 600-cell is the finite regular four-dimensional polytope composed of 600 tetraheda. Its skeleton $G=(V,E)$ has $|V|=120$ vertices and $|E| = 720$ edges. \textsc{Max-Cut}(G) seems to be unknown. The Goemans-Williamson algorithm run over many instances yields $\textsc{Max-Cut}(G) \geq 432$. BUR02 improves this to $\textsc{Max-Cut}(G) \geq 436$. Our approach does not further improve on this and also shows $\textsc{Max-Cut}(G) \geq 436$ (though the expected size of the cut increases). However, looking at the final configuration of points (see Fig. 6), the final configurations are quite different.

\vspace{-10pt}

 \begin{center}
\begin{figure}[h!]
\begin{tikzpicture}
\node at (-1,0) {\includegraphics[width=0.30\textwidth]{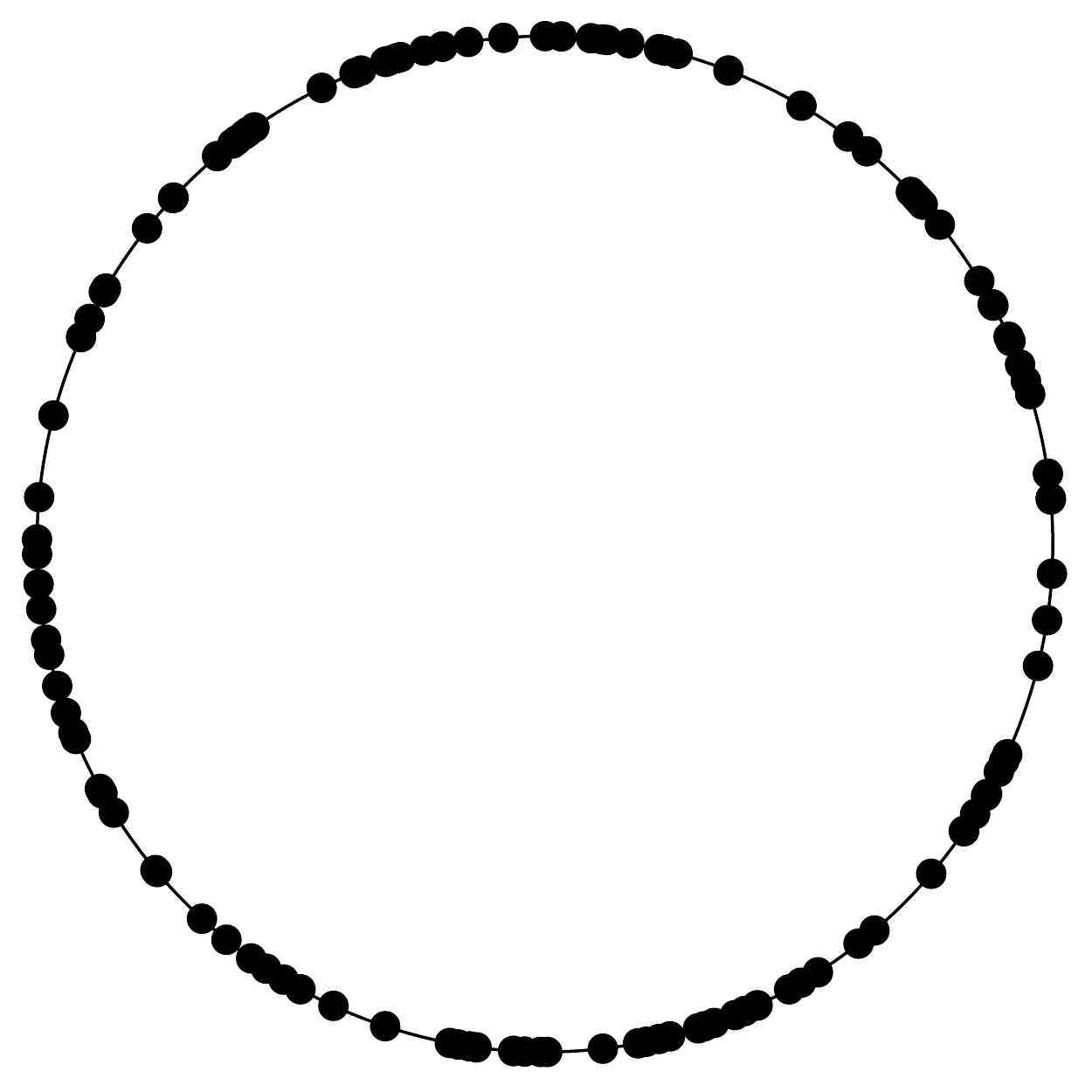}};
\node at (5,0) {\includegraphics[width=0.30\textwidth]{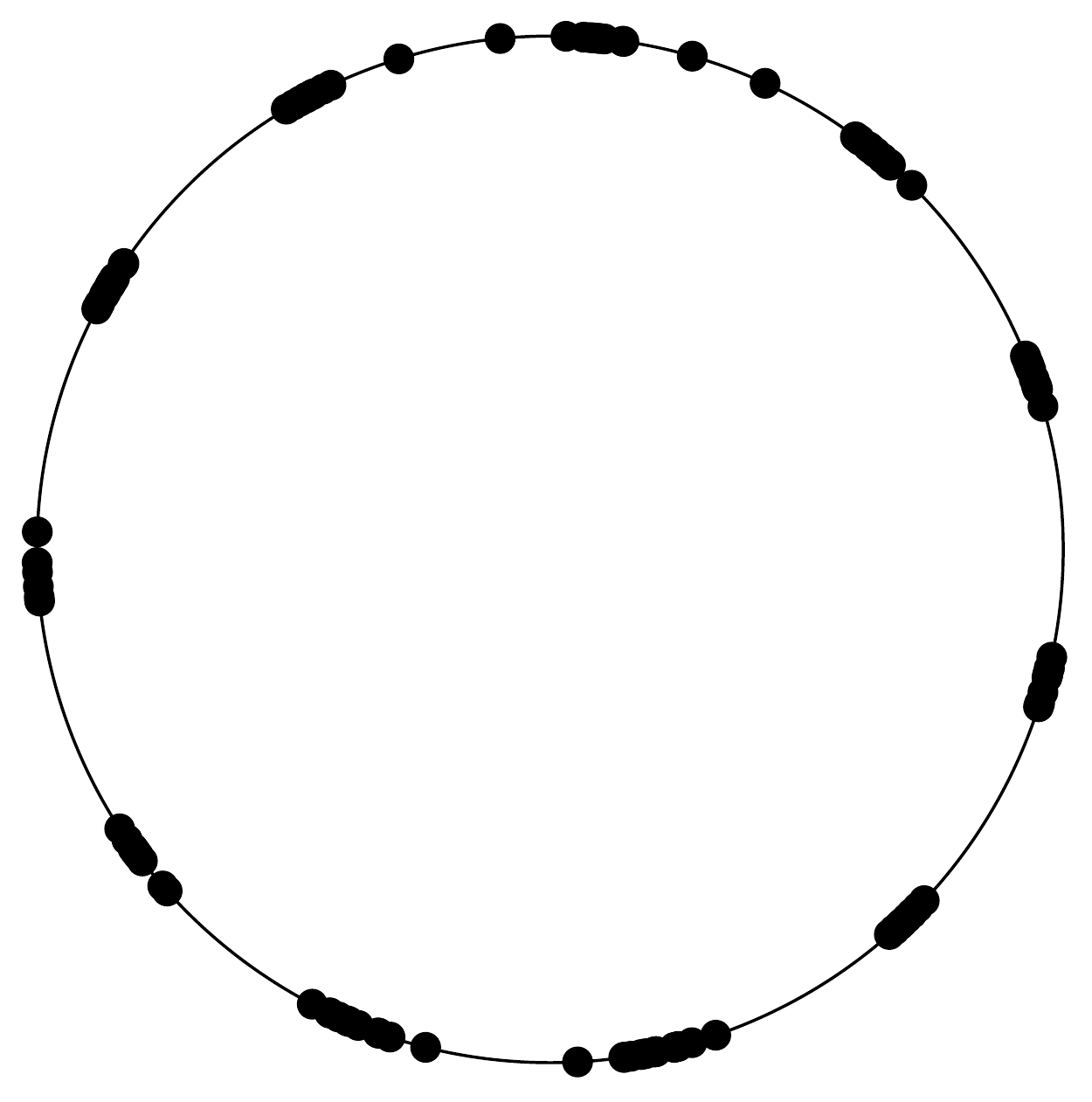}};
\end{tikzpicture}
\caption{Final configuration of BUR02 (left) and our method (right) when applied to the skeleton graph of the 600-cell. Both configurations show $\textsc{Max-Cut}(G) \geq 436$.}
\end{figure}
\label{fig:mont}
\end{center}

\vspace{-10pt}

We believe the reason for this `crystallization' phenomenon to be the following: for the particular function $g$ that we use, we have
$$ \mathbb{E}_{} ~\mbox{edges} \geq 0.973 \cdot \textsc{Max-Cut}(G)$$
for the minimal energy configuration (though, of course, we cannot be sure of having found a minimal energy configuration). This has a very powerful implication because it means the \textit{for virtually every line}, the induced partition is necessarily \textit{very close} to \textsc{Max-Cut}.  Whenever the arising distribution of points is not simply concentrated at two antipodal points, then the final configuration has to show different ways how a partition of vertices with a number of edges close to \textsc{Max-Cut} can be achieved. We believe that this explains the arising crystallization that we observe. It also indicates that this such minima should actually induce a rather interesting ordering of the vertices of the Graph in terms of groups that have strong interactions with antipodal groups. One would assume that the highly structured picture in Fig. 6 somehow reflects the underlying structure of the Graph.

\subsection{Which $g$ should one use?} One important question is the choice of the function $g$. Our main result suggests that we should pick 
$$ g(x) \sim 1 - \frac{2}{\pi} d_{\mathbb{S}^1}(0,x)$$
so that global minima correspond to a good approximation of \textsc{Max-Cut}. However, this is counter-balanced by optimization concerns -- the global minimum having good properties will not be of any use to us if we cannot get close to it. At this point, we have no good theoretical reason to choose any particular $g(x)$ and we believe this to be an interesting problem.
\begin{quote}
\textbf{Question.} What are good choices for $g$?  Which properties of $g$ lead to the functional having `nice' energy landscapes?
\end{quote}

We found that $g(x)$ being close to $1 - (2/\pi) \cdot d_{\mathbb{S}^1}(x,0)$ is indeed beneficial for the quality of the solution but also makes optimization harder. Smooth functions tend to be easier to optimize, hence our choice to use a truncated Fourier series approximation of $1 - (2/\pi) \cdot d_{\mathbb{S}^1}(x,0)$. In practice, it may well be a good idea to initialize with some $g$ and then change the choice of $g$ after a while. \\
\textit{Trigonometric Polynomials.} We mention one particular reason that might speak in favor if using trigonometric polynomials (and is completely unconnected to any considerations about the energy landscape). Fix $\left\{ \theta_1, \dots, \theta_n \right\} \subset \mathbb{S}^1$.  We can pick an arbitrary $\theta_i$, keep the remaining angles fixed and ask ourselves how the function
$$  \sum_{i,j=1}^{n} a_{ij} \cdot g(\theta_i - \theta_j) \qquad \mbox{behaves as a function of}~\theta_i.$$
If $g$ is a trigonometric polynomial of degree $d$, then this sum is, as a function of $\theta_i$, also a trigonometric polynomial of degree $d$ because trigonometric polynomials of degree $d$ are an invariant subspace under translation. This means that this function, as a function $\theta_i$, is globally quite simple and we can find its global minimum. This is particularly striking in the case of BUR02: the function
$$ h(\theta_i) =  \sum_{j=1 \atop j \neq i}^{n} a_{ij} \cdot \cos(\theta_i - \theta_j) +  \sum_{j=1 \atop j \neq i}^{n} a_{ji} \cdot \cos(\theta_j - \theta_i)$$
is a function of the form
$$ h(\theta_i) = A \cdot \cos{(\theta_i -B)},$$
where $A$ and $B$ depend on all the other variables. However, such a function is very easy to minimize globally: set $\theta_i = B + 3 \pi/2 ~(\mbox{mod}~2\pi$).  This persists when passing from the cosine to trigonometric polynomials of degree $d$ (which is a $2d-$dimensional vector space with rather nicely behaved functions in it that always have a lot of structure and are easier to minimize than generic functions). This allows for non-local optimization schemes along the following lines: pick a variable $\theta_i$, freeze all the other varables, compute where one would place $\theta_i$ to minimize the energy and move it there. One would expect that the effectiveness of such a scheme depends on the function $g(x)$ which brings us back to the question raised above.

\section{Proof}

\begin{proof} We will now prove the Theorem. The argument is identical to the classical randomized rounding argument except that we are working with an arbitrary function $g$ and track its dependence.
Suppose the Max-Cut solution is given by the splitting $V = A \cup B$. Then we can set all the vertices in $A$ to have angle $\theta_a = 0$ and all the angles in $B$ to have $\theta_b = \pi$ and compute the energy of this configuration. There are $\textsc{Max-Cut}(G)$ edges getting weight $-1$ and $|E| - \textsc{Max-Cut}(G)$ edges getting weight 1. Every edge is counted twice, therefore
\begin{align*}
\min_{\theta_1, \dots, \theta_n} f(\theta_1, \dots, \theta_n) \leq 2 \cdot |E| - 4 \cdot \textsc{Max-Cut}(G).
 \end{align*}
Suppose conversely that we have a configuration with small energy given by the configuration of angles $\left\{\theta_1, \theta_2, \dots, \theta_n \right\} \subset [0,2\pi]$. The likelihood of two specific vertices $i,j \in V$ ending up in different partitions is given by the likelihood of $\theta_i$ and $\theta_j$ being cut by a hyperplane. That quantity has a simple expression given by
$$ \mathbb{P}\left( \theta_i, \theta_j~\mbox{in different halfspaces}\right) = \frac{|\theta_i - \theta_j|_{\mathbb{S}^1}}{\pi},$$
 where $| \cdot |$ denotes the shortest distance on $\mathbb{S}^1$ (and is thus always less than $\pi$). Using linearity of expectation, we can compute
 the expected number of edges across a randomly chosen line
  \begin{align*}
\mathbb{E}~\mbox{edges} &=   \frac{1}{2} \sum_{i,j=1}^{n} a_{ij} \cdot \mathbb{P}\left( \theta_i, \theta_j~\mbox{in different halfspaces}\right) \\
&=  \frac{1}{2} \sum_{i,j=1}^{n} a_{ij} \cdot \frac{|\theta_i - \theta_j|_{\mathbb{S}^1}}{\pi}.
   \end{align*}
 At this point, we use that the distance function satisfies, tautologically,
 $$ \frac{|\theta_i - \theta_j|_{\mathbb{S}^1}}{\pi} \geq \left(  \min_{0 \leq x \leq \pi}   \frac{2}{\pi}  \frac{x}{1 -g(x)}  \right) \cdot \frac{1-g(\theta_i - \theta_j)}{2},$$
 we have
 $$ \frac{1}{2} \sum_{i,j=1}^{n} a_{ij}  \frac{|\theta_i - \theta_j|}{\pi} \geq \left(  \min_{0 \leq x \leq \pi}   \frac{2}{\pi}  \frac{x}{1 -g(x)}  \right) \sum_{i,j=1}^{n} a_{ij} \frac{1- g(\theta_i- \theta_j)}{4}.$$
 This sum simplifies to
 \begin{align*}
  \sum_{i,j=1}^{n} a_{ij} \frac{1- g(\theta_i- \theta_j)}{4} &= \frac{|E|}{2} - \frac{1}{4} \sum_{i,j = 1}^{n} a_{ij} g(\theta_i - \theta_j) \\
 &=  \frac{|E|}{2} - \frac{1}{4} f(\theta_1, \dots, \theta_n).
  \end{align*}
Therefore, we have
 $$ \mathbb{E}~\mbox{edges} \geq  \left(  \min_{0 \leq x \leq \pi}   \frac{2}{\pi}  \frac{x}{1 -g(x)}  \right) \left( \frac{|E|}{2} - \frac{1}{4} \sum_{i,j = 1}^{n} a_{ij} \cos{(\theta_i - \theta_j)} \right).$$ 
 The remaining question is simply how small we can make this Kuramoto-type energy: by the argument above, we have
$$\min_{\theta_1, \dots, \theta_n} f(\theta_1, \dots, \theta_n) \leq 2 \cdot |E| - 4 \cdot \textsc{Max-Cut}(G).$$
Thus, if $\left\{\theta_1, \dots, \theta_n\right\}$ is a minimal energy configuration of the Kuramoto energy, 
 $$ \mathbb{E}~\mbox{edges} \geq  \left(  \min_{0 \leq x \leq \pi}   \frac{2}{\pi}  \frac{x}{1 -g(x)}  \right)\cdot \textsc{Max-Cut}(G)$$
 which completes the argument. 
\end{proof}

\end{document}